\documentclass{article}
\usepackage{amsmath}
\usepackage{amssymb}
\usepackage{color}

\usepackage{latexsym}
\usepackage[english]{babel}
\usepackage{graphicx}
\usepackage{verbatim}
\usepackage[OT1]{fontenc} 
\usepackage[latin1]{inputenc}

\newcommand{\llabel}[1]{{\label{#1}}}
\newcommand{\ex}[1]{}

\font\tenmsb=msbm10 \font\sevenmsb=msbm7 \font\fivemsb=msbm5
\newfam\msbfam
\textfont\msbfam=\tenmsb \scriptfont\msbfam=\sevenmsb
\scriptscriptfont\msbfam=\fivemsb
\def\Bbb#1{{\fam\msbfam\relax#1}}

\newcommand{\bi}{\begin{itemize}}
\newcommand{\ei}{\end{itemize}}
\newcommand{\bd}{\begin{description}}
\newcommand{\ed}{\end{description}}

\newcommand{\bqn}{\begin{eqnarray}}
\newcommand{\eqn}{\end{eqnarray}}
\newcommand{\eqnn}{\nonumber\end{eqnarray}}
\newcommand{\eqnl}[1]{\llabel{#1}\end{eqnarray}}

\newcommand{\nn}{\nonumber}
\newcommand{\ba}[1]{\begin{array}{#1}}
\newcommand{\ea}{\end{array}}

\newcommand{\R}{\Bbb{R}}

\newcommand{\fine}{\end{document}}

\def \trait (#1) (#2) (#3){\vrule width #1pt height #2pt depth #3pt}
\def \qed{\hfill
        \trait (0.1) (6) (0)
        \trait (6) (0.1) (0)
        \kern-6pt
        \trait (6) (6) (-5.9)
        \trait (0.1) (6) (0)
\medskip}
\def \qedmio{\hfill
             \trait (8) (8) (-0.1)
             \medskip}

\newtheorem{ml}{\bf Lemma}
\newtheorem{Theorem}{\bf Theorem}
\newtheorem{mo}{\bf \underline{{\sl Observation}}}
\newtheorem{mcc}{\bf Corollary}
\newtheorem{Definition}{\bf Definition}
\newtheorem{mpr}{\bf Proposition}
\newtheorem{mproperty}{\bf Property}

\newcommand{\bt}{\begin{Theorem}}
\newcommand{\et}{\end{Theorem}}
\newcommand{\bl}{\begin{ml}}
\newcommand{\el}{\end{ml}}
\newcommand{\bo}{\noindent\begin{mo}\rm}
\newcommand{\eo}{\end{mo}}
\newcommand{\bp}{\begin{mpr}}
\newcommand{\ep}{\end{mpr}}
\newcommand{\bc}{\begin{mcc}}
\newcommand{\ec}{\end{mcc}}
\newcommand{\bdeff}{\begin{Definition}}
\newcommand{\edeff}{\end{Definition}}
\newcommand{\bproperty}{\begin{mproperty}}
\newcommand{\eproperty}{\end{mproperty}}

\newtheorem{mrem}{\bf \underline{{\sl Remark}}}
\newcommand{\brem}{\begin{mrem}\rm}
\newcommand{\erem}{\end{mrem}}

\newcommand{\lam}{\lambda}
\newcommand{\g}{\gamma}
\newcommand{\al}{\alpha}
\newcommand{\eps}{\varepsilon}
\newcommand{\de}{\delta}

\newcommand{\F}{{\cal F}}

\newcommand{\Z}{{\mathcal Z}}

\newcommand{\U}{{\mathcal U}}

\newcommand{\ga}{\gamma}

\newcommand{\ph}{\varphi}
\renewcommand{\c}{{\cos(\ph)}}

\newcommand{\guas}{{\bf GUAS}}
\newcommand{\lf}{{\bf LF}}

\newcommand{\two}{two-dimensional bilinear switched system}

\newcommand{\lunga}{\longrightarrow}

\newcommand{\del}{\delta}

\newcommand{\tr}{{\rm tr}}

\newcommand{\para}{\frac{2}{\sqrt{|\delta_{A_1}|}   }}
\newcommand{\parb}{\frac{2}{\sqrt{|\delta_{A_2}|}   }}

\begin{document}
\begin{center} \noindent
{\LARGE{\sl{\bf  
A note on stability conditions for planar switched systems}}}\footnote{The first two authors were supported by a FABER grant of Universit\'e de Bourgogne}
\end{center}

\vskip 1cm
\begin{center}
Moussa BALDE,

{\footnotesize LMDAN-LGDA D\'epartement de Math\'ematiques et Informatique,
UCAD, Dakar-Fann, Senegal\\
{mbalde@ucad.sn}}
\end{center}
\begin{center}

\vskip .2cm Ugo BOSCAIN,

{\footnotesize Le2i, CNRS, Universit\'e de Bourgogne, B.P. 47870,
21078
Dijon Cedex, France\\ ugo.boscain@u-bourgogne.fr}

\end{center}
\begin{center}

\vskip .2cm  Paolo MASON

{\footnotesize Istituto per le Applicazioni
del Calcolo Mauro Picone - CNR Viale del Policlinico 137 - 00161
Roma, Italy\\
and\\
Laboratoire des signaux et syst\`emes, Universit\'e Paris-Sud, CNRS, Sup\'elec, 91192 Gif- 
Sur-Yvette, France\\ Paolo.Mason@lss.supelec.fr
}
\end{center}

 \medskip

\begin{abstract}
This paper is concerned with the stability problem for the planar linear switched system $\dot
x(t)=u(t)A_1x(t)+(1-u(t))A_2x(t)$, where the real matrices $A_1,A_2\in \R^{2\times 2}$ are Hurwitz
and $u(\cdot):[0,\infty[\to\{0,1\}$ is a measurable function. We give  coordinate-invariant necessary and sufficient conditions on $A_1$ and $A_2$ under which the system is asymptotically
stable for arbitrary switching functions $u(\cdot)$.  The new conditions unify those given in previous papers and are simpler to be verified since we are reduced to study 4 cases instead of 20.
Most of the cases are analyzed in terms of the function $\Gamma(A_1,A_2)=\frac{1}{2}(\tr(A_1)
\tr(A_2)- \tr(A_1A_2))$.

\end{abstract}

\vspace{10pt}


\section{Introduction}

\vspace{10pt}

Let $A_1$ and $A_2$ be  two $2\times 2$ real Hurwitz matrices. In this paper we are concerned with the problem of finding necessary and sufficient conditions on $A_1$ and $A_2$ under which the switched system
\bqn \dot
x(t)=u(t)A_1x(t)+(1-u(t))A_2x(t), ~~x=(x_1,x_2)\in\R^2,\eqnl{sw-ls}
is  globally asymptotically stable, uniformly with respect to measurable switching functions  $u(\cdot):[0,\infty[\to\{0,1\}$ ({\bf GUAS} for short, see Definition~\ref{d-stability} below).

This problem has been studied in \cite{sw-1}  in  the case in which both $A_1$ and $A_2$ are diagonalizable in
$\Bbb{C}$ ({\it diagonalizable} case in the following) and in \cite{sw-balde} in the case in which at least one among  $A_1$ and  $A_2$ is not ({\it nondiagonalizable} case in the following). (See also \cite{sw-lyapunov} as well as the related work \cite{marga}.)

In both cases the stability conditions are given in terms of coordinate-invariant parameters. Unfortunately the parameters used in the diagonalizable case become singular in the nondiagonalizable one and therefore the two cases were studied separately.

The purpose of this note is to unify and simplify these conditions, reformulating  them in terms of new invariants that permit to treat all cases  at the same time.

We have reduced the cases to be studied from 20
(14 in the diagonalizable case\footnote{The stability conditions given in \cite{sw-1} were not  correct in the case called {\bf RC.2.2.B}.
See \cite{sw-lyapunov} for the correction}  and 6 in the nondiagonalizable one) to the following 4 cases (see Theorem~\ref{t-main}).

\begin{description}

\item[S1:] the first one corresponds to the case in which there exists a common quadratic Lyapunov function. The condition of {\bf S1} is indeed equivalent to the condition given in \cite{shorten} but is simpler to check. Recall however that the existence of a common quadratic Lyapunov function is only a sufficient condition for \guas\ (i.e.  there exist \guas\ systems not admitting a quadratic Lyapunov function). See \cite{DM, sw-lyapunov} for details.

\item[S2:] the second  one corresponds to the situation in which there exists
$v\in(0,1)$ such that $vA_1+(1-v)A_2$ has a positive real eigenvalue.
In this case the system is unbounded since it is possible to build a trajectory going to infinity approximating the (non admissible) trajectory  corresponding to $u(t)\equiv v$ and having the direction of the unstable eigenvector of $vA_1+(1-v)A_2$.

\item[S3:] in the third case there exists a nonstrict common quadratic Lyapunov function. The system is only uniformly stable, but not {\bf GUAS}, since there exists a trajectory not tending to the origin when $t$ goes to infinity.

\item[S4:] in the fourth case the stability analysis of the system reduces to the study of a single trajectory called {\it worst trajectory}. If this trajectory tends to the origin then the system is \guas\ (in this case there  exists a polynomial Lyapunov function, but not a quadratic one). If it is periodic then the system is uniformly stable but not \guas. If it is unbounded then the system is unbounded.

\end{description}

For a discussion of various issues related to stability of switched systems, we
refer the reader to \cite{DM,liberzon-book}. 

The paper is organized as follows.
In Section~\ref{s-stab} we recall the fundamental notions of stability and the different types of Lyapunov functions used in the paper. Section~\ref{s-2} contains our main result. In Section~\ref{s-normal} we define the normal forms that are needed in the proof. In Section~\ref{s-proof}  we give the details of the proof.

\subsection{Notions of stability}
\label{s-stab}
Let us recall some classical notions of stability which will be used in the
following.
\bdeff \label{d-stability} For $\de>0$ let $B_\de$ be the
unit ball of radius $\de$, centered in the origin. Denote by $\U$
the set of measurable functions defined on $[0,\infty[$  and taking
values on $\{0,1\}$. Given $x_0\in \R^2$, we denote by $\g_{x_0,u(\cdot)}(\cdot)$
the trajectory of \eqref{sw-ls} based in $x_0$ and corresponding to the
control $u(\cdot)$.
We say that the system \eqref{sw-ls} is 
\bi
\item {\bf unbounded} at the origin if there exist $x_0\in \R^2$
and $u(\cdot)\in\U$ such that $\ga_{x_0,u(\cdot)}(t)$ goes to infinity as $t$ goes to infinity;
\item {\bf uniformly stable} at the origin if for every $\eps>0$ there exists
$\de>0$ such that  $\g_{x_0,u(\cdot)}(t)\in B_{\eps}$ for every $t>0$, for every $u(\cdot)\in\U$ and every $x_0\in B_{\delta}$;
\item {\bf globally uniformly asymptotically stable} at the origin ({\bf
GUAS}, for short) if it is uniformly stable at the origin and globally uniformly
attractive,  i.e., for every $\de_1, \de_2>0$, there exists $T>0$
such that $\ga_{x_0,u(\cdot)}(t)\in B_{\de_1}$ for every $t\geq T$, for every $u(\cdot)\in\U$
and every $x_0\in B_{\de_2}$.
\ei
\edeff
\brem
\label{r-conv}
The stability properties of the system  \eqref{sw-ls} do not change if we allow measurable switching functions taking values in $[0,1]$ instead of $\{0,1\}$ (see for instance \cite{sw-lyapunov}).
More precisely
the system \eqref{sw-ls} with $u(\cdot):[0,\infty[\to \{0,1\}$ is \guas\
(resp. uniformly stable, resp. unbounded) if and only the system
\eqref{sw-ls} with $u(\cdot):[0,\infty[\to [0,1]$ is.
In the following we name {\it convexified system} the switched system with $u(\cdot)$
taking values in $[0,1]$.
\erem

Since the stability properties of the system \eqref{sw-ls} do not depend on the
parametrization of the integral curves of $A_1x$ and $A_2x$, we have the following.
\bl
\label{stin} If the switched system $\dot
x=u(t)A_1x+(1-u(t))A_2 x,~~~u(\cdot):[0,\infty[\to\{0,1\},$ has one of the
stability properties given in Definition \ref{d-stability}, then the
same stability property holds for the system $\dot
x=u(t)(A_1/\al_1)x+(1-u(t))(A_2/\al_2)x$, for every $\al_1,\al_2>0$.
\el

\bdeff \llabel{d-CLF} A \emph{common Lyapunov function}
(\ {\bf LF} for short) $V:\R^2\lunga \R^+$ for a switched
system of the form \eqref{sw-ls} is a continuous function such that
$V(\cdot)$ is positive definite (i.e. $V(x)>0$, $\forall x\neq 0$,
$V(0)=0$) and $V(\cdot)$ is strictly decreasing along nonconstant
trajectories.

A positive definite continuous function $V:\R^2\lunga \R^+$ is said to be a \emph{nonstrict common Lyapunov function} if $V(\cdot)$ is nonincreasing along nonconstant trajectories.

A  \emph{common quadratic Lyapunov function} (quadratic \lf\ for short) is a function of the form $V(x)=x^T Px$ where $P$ is a positive definite symmetric matrix and the matrices $A_1^T P+PA_1$ and  $A_2^T P+PA_2$ are negative definite.
\edeff

We recall that, for systems of type \eqref{sw-ls},  the existence of a  \lf\  is equivalent to \guas \footnote{In \cite{blanchini,sw-lyapunov,molch1,molch2} it is actually shown that the \guas\ property is  equivalent to the existence of a polynomial \lf} (see for instance \cite{DM}). 
Moreover the existence of a nonstrict \lf\ guarantees the uniform stability of \eqref{sw-ls}.

\section{Stability conditions for  \two s}
\label{s-2}
We start this section  by defining the notations and the objects that will be used to state our stability result. In the following the word invariant will indicate any object which is invariant with respect to coordinate transformations.
As usual, we denote by $\det(X)$ and $\tr(X)$ the determinant and the trace of a matrix $X$. If $X\in \R^{2\times 2}$ the discriminant is defined as
$$\delta_X =\tr(X)^2-4\, \det(X)\,.$$
Given a pair of matrices $X,Y$ we define the following object:
$$\Gamma(X,Y):=\displaystyle\frac{1}{2}(\tr(X) \tr(Y)- \tr(X Y))\,.$$
By means of these invariants we can define the following invariants associated to \eqref{sw-ls}:

\bqn
& \tau_i:=\left\{\ba{l}\displaystyle\frac{\tr(A_i)}{\sqrt{|\delta_{A_i}|}}\quad\mbox{ if }\ \del_{A_1}\neq 0\,,\ \del_{A_2}\neq 0\,,\\
\displaystyle\frac{\tr(A_i)}{\sqrt{|\del_{A_j} |}}\quad\mbox{ if }\ \del_{A_1} \del_{A_2} = 0 \mbox{ but } \del_{A_j} \neq 0\,, \\
\displaystyle\frac{\tr(A_i)}{2}\quad\mbox{ if }\ \del_{A_1}=\del_{A_2}=0\,,\ea\right.& \nn\\
& k := \displaystyle\frac{2\tau_1\tau_2}{\tr(A_1)\tr(A_2)}\Big(\tr(A_1 A_2)-\frac12\tr(A_1) \tr(A_2)\Big)\,,& \nn\\
&\Delta :=
4(\Gamma(A_1,A_2)^2-\Gamma(A_1,A_1)\Gamma(A_2,A_2))\,.& \nn \\
&\mathcal{R}:=\displaystyle\frac{2\Gamma(A_1,A_2)+\sqrt{\Delta}}{2\sqrt{\det(A_1)\det(A_2)}} e^{\tau_1 t_1+\tau_2 t_2}\,, \mbox{where, for }\ i=1,2\,,& \nn\\
&t_i:=\left\{ \ba{lr} \frac{\pi}2-\mbox{arctan}\frac{\tr(A_1)\tr(A_2)(k\tau_i+\tau_{3-i})}{2 \tau_1\tau_2\sqrt{\Delta}} & \mbox{ if }\ \delta_{A_i}<0 \\
\mbox{arctanh}\frac{2 \tau_1\tau_2\sqrt{\Delta}}{\tr(A_1)\tr(A_2)(k\tau_i-\tau_{3-i})} & \mbox{ if }\ \delta_{A_i}>0\\
\frac{2\sqrt{\Delta}}{(\tr(A_1 A_2)-\tr(A_1)\tr(A_2)/2)\tau_i} &  \mbox{if }\ \delta_{A_i}=0\ea\right.\,.& \nn
\eqn

\brem 
Let us define
$$
sign(x):=\left\{\ba{r}+1\ \mbox{ if $x>0$  }\\
0\ \mbox{ if $x=0$ }\\
-1\ \mbox{ if $x<0$
  } \ea\right.\,.
$$
Notice that, for every matrix $X\in\R^{2\times 2}$, one has $\Gamma(X,X)=\det (X)$. Also, since the Killing form of $\R^{2\times 2}$ is defined as $\mathcal{B}(X,Y)=4\,\tr (XY)-2\,\tr(X)\tr(Y)$ one has $sign(k)=sign(\mathcal{B}(A_1,A_2))$. Finally, notice that if $A_1,A_2$ are Hurwitz then $\tau_i<0\,$ for $\,i=1,2\,$ and  $\ sign\, \Gamma(A_1,A_2)=sign(\tau_1\tau_2-k)$.
\erem

\subsection{Statement of the results}
\label{s-main}
In this section we state our main result which characterizes conpletely the stability properties of \two s.
Our necessary and sufficient conditions apply both to the non-degenarate cases studied in \cite{sw-1} and to the degenarate ones studied in \cite{sw-balde}.

\begin{Theorem}\label{t-main}
We have the following stability conditions for the system \eqref{sw-ls}
 \bd
     \item[S1]  If \  $\Gamma{(A_1,A_2)} >  -\sqrt{\det(A_1)\det(A_2)}, \  \text{and} \ \tr(A_1 A_2)>  -2\sqrt{\det(A_1)\det(A_2)}$ then the system admits a quadratic \lf.

    \hspace{-.35 cm} If  ~~$-\sqrt{\det(A_1)\det(A_2)}<\Gamma(A_1,A_2)\leq\sqrt{\det(A_1)\det(A_2)}$ then the condition $\tr(A_1 A_2)>  -2\sqrt{\det(A_1)\det(A_2)}$ is automatically satisfied. As a consequence the system admits a quadratic \lf.
   \item[S2] If \ $\Gamma{(A_1,A_2)} < -\sqrt{\det(A_1)\det(A_2)},$ then the
     system is  unbounded,
    \item[S3] If \ $\Gamma{(A_1,A_2)}= - \sqrt{\det(A_1)\det(A_2)},$ then the system is uniformly stable but not \guas,
     \item[S4] If \  $\Gamma{(A_1,A_2)} >  \sqrt{\det(A_1)\det(A_2)},\ \ \text{and} \ \ \tr(A_1 A_2)\leq  -2\sqrt{\det(A_1)\det(A_2)}$ then the system is \guas,
uniformly stable (but not \guas) or unbounded respectively if
                $$ {\mathcal R} <1,{\mathcal R} =1, {\mathcal R} >1 .$$

      \ed
\end{Theorem}
The following corollary will be derived from item {\bf S1} of the previous theorem.
\bc
If $\det([A_1,A_2])\geq 0$ then the system admits a quadratic \lf.
\label{pos}
\ec

\brem 
In the diagonalizable case $\del_{A_1} \del_{A_2}\neq 0$ the parameters $\tau_1$, $\tau_2$, and $k$ are invariant under the transformation $(A_1,A_2) \to(A_1/\al_1 , A_2/\al_2)$,
for every $\al_1,\al_2>0.$ This is no more true in the nondiagonalizable case.
Notice however that in any case the stability conditions of Theorem~\ref{t-main} do not depend on coordinate transformations or on rescalings of the type $(A_1,A_2) \to(A_1/\al_1 , A_2/\al_2)$. This is true in particular for the function $\mathcal{R}$.
\label{rescal}
\erem

\section{Proof of the main results}

\subsection{Normal forms}
\label{s-normal}
The aim of this section is to reduce all the possible choices of the matrices $A_1,A_2$ to suitable normal forms, obtained up to coordinates transformations and rescaling of the matrices (see Lemma~\ref{stin} and Remark~\ref{rescal} above), and depending directly on the coordinate invariant parameters introduced above. The normal forms used  here describe all the possible situations for two-dimensional bilinear switched systems, covering at the same time the diagonalizable case studied in \cite{sw-1} and the nondiagonalizable one studied in \cite{sw-balde}.
They will play a key role in the proof of our results.
\bl
We have the following cases depending on the rank of $[A_1,A_2]$:
\begin{enumerate}
\item If $\det([A_1,A_2])\neq 0$, up to a linear change of coordinates and a renormalization according to Lemma~\ref{stin}, we can assume the following.
\bqn &&A_1=\left(\ba{cc}\tau_1&1\\
sign(\delta_{A_1})&\tau_1\ea\right),\label{e-nA}\eqn
\begin{enumerate}
  \item If $\det([A_1,A_2]) < 0$ there exists $F \in \mathbb{R}$, $|F|\geq 1$ 
  such that
  $$ F+ \frac{sign(\delta_{A_1} \delta_{A_2})}{F}=2k$$ and $A_2$ has the form
  \bqn
&&A_2=\left(\ba{cc}\tau_2&sign(\delta_{A_2})/F\\F&\tau_2\ea\right)\label{e-nB1}.
\eqn
  \item If  $\det([A_1,A_2]) > 0$ then $\delta_{A_i}>0$ for $i=1,2$,  $k \in (-1,1)$ and $A_2$ has the form
\bqn  &&A_2=\left(\ba{cc}\tau_2+
\sqrt{1-k{^2}}&k\\k&\tau_2-\sqrt{1-k{^2}}\ea\right)\label{e-nB2}.
\eqn
\end{enumerate}
\item Under the hypothesis that
$rank([A_1,A_2])=1$, it is always possible, up to exchanging $A_1$ and $A_2$,
to find a linear change of coordinates which diagonalizes $A_1$ and
renders $A_2$  upper triangular.
\item If $[A_1,A_2]=0$ then it must be $sign(\delta_{A_1})=sign(\delta_{A_2})$. If $\delta_{A_i}<0$ for $i=1,2$ or $\delta_{A_i}>0$ for $i=1,2$ then,  up to a linear change of coordinates and a renormalization according to Lemma~\ref{stin}, $A_1,A_2$ assume the form \eqref{e-nA} and \eqref{e-nB1}, respectively, with $F=k=\pm 1$.  If $\delta_{A_1}=\delta_{A_2}=0$  then $A_1,A_2$ can be put in upper triangular form with the elements of $A_i$ equal to $\tau_i$ on the diagonal, for $i=1,2$.
\end{enumerate}
\label{l-nf-r}
\el

{\bf Proof of Lemma \ref{l-nf-r}.}
For simplicity we will prove the lemma just in the case $\det([A_1,A_2])\neq 0$, the other case being analogous.
Note that Lemma \ref{l-nf-r} was proven in \cite{sw-1} when $\delta_{A_1}<0,\  \delta_{A_2}<0$ and in \cite{sw-balde} in the case $\delta_{A_1} \delta_{A_2}=0$. Therefore, we can assume either $\delta_{A_1}>0$ or $\delta_{A_2}>0$. First consider the case $\delta_{A_1}>0$.  In  this case we can find a system of coordinates such that
\bqn &&A_1=\left(\ba{cc}\lam_1&
0\\0&\lam_2\ea\right),~~~A_2=\left(\ba{cc}a&b\\c&d\ea\right),
~~a,b,c,d\in\R. \eqn Without any loss of generality we
 can assume that $\lam_1< \lam_2.$ The discriminant of $A_1$ is
 $(\lam_2-\lam_1)^2$ and the discriminant of $A_2$ is
$\delta_{A_2}=(a-d)^2+4 bc$, which can be positive or negative. We have
\bqn [A_1,A_2]=\left(\ba{cc}0&{b(\lam_1 - \lam_2)}\\c(\lam_2 -
\lam_1)&0\ea\right). \eqn
 Therefore  $\det{[A_1,A_2]} = bc (\lam_1 -\lam_2)^2 \ \ \text{and}
  \ \  \delta_{[A_1,A_2]}= -4bc(\lam_1 -\lam_2)^2$  
  and $k= \frac{d-a}{\sqrt{|\delta_{A_2}|}}$.
  If $\det{[A_1,A_2]}<0$ consider the linear transformation $$T =\left(\ba{cc} -\sqrt{\frac{-b}{c}}&
  \sqrt{\frac{-b}{c}}\\1&1\ea\right)\,,$$ which diagonalizes $[A_1,A_2]$. Then a straightforward computation shows that
$$
\para T^{-1}A_1T  =  \left(\ba{cc}\tau_1&1\\1&\tau_1 \ea\right) \qquad \mbox{and}  \qquad \frac{2}{\sqrt{|\delta_{A_2}|}} T^{-1}A_2T  =  \left(\ba{cc} \tau_2 & sign(\delta_{A_2})/F\\ F&\tau_2\ea\right),
$$
where $F$ satisfies the equation  $F+sign(\delta_{A_2})/F=2k$,
 and moreover we can assume $|F|\geq 1$ up to eventually exchange the reference coordinates.
If $\delta_{A_1}<0$ then it must be $\delta_{A_2}>0$ and, exchanging the roles of $A_1$ and $A_2$, we can repeat the previous procedure obtaining
$$A_1=\left(\ba{cc} \tau_1 & sign(\delta_{A_1})/F\\ F&\tau_1\ea\right)\,,~~~~A_2=\left(\ba{cc} \tau_2 & 1 \\ sign(\delta_{A_2}) &\tau_2\ea\right).$$
Then the required normal forms are obtained by exchanging the coordinates and by a dilation along one of the coordinate axis.

Consider now the case $\det{[A_1,A_2]}= bc(\lam_1 -\lam_2)^2 >0.$
We have $$b c>0~~\Rightarrow~~\delta_{A_2}=(a-d)^2+4 bc>0  \ \  \text{and} \ \
|k|=\frac{|a-d|}{\sqrt{(a-d)^2+ 4 bc}}<1\,.$$ In this case $[A_1, A_2]$ is
no more diagonalizable. Using the transformation $$U =\left(\ba{cc} \frac1c \sqrt{bc} & -\frac1c \sqrt{bc} \\ 1&1\ea\right)$$
we get
$$ \para U^{-1} A_1 U =  \left(\ba{cc} \tau_1 & 1 \\ 1 &\tau_1 \ea\right) \qquad \mbox{and} \qquad \parb U^{-1} A_2 U = \left(\ba{cc} \tau_2+\sqrt{1-k^2} & k \\ k&\tau_2-\sqrt{1-k^2}\ea\right),$$
which concludes the proof of the lemma.\qedmio


\subsection{Proof of Theorem~\ref{t-main} and Corollary~\ref{pos}}
\label{s-proof}
To prove our main result we will assume, from now on, that $A_1,A_2$ are under the normal forms given by Lemma~\ref{l-nf-r}.
The following lemma, which can be proved by direct computation, will be used to take advantage of the conditions of \cite{shorten} which describe the systems admitting a quadratic \lf.
\bl
 For any $\sigma \in [0,1],$ we define $$\phi(\sigma):=\det{(\sigma A_1 + (1-\sigma) A_2)}\,,  \qquad \psi(\sigma):= \det{(\sigma A_1 + (1-\sigma) A_2^{-1})}\,. $$
We have
\bqn
\phi(\sigma)= \sigma^2 \det{A_1}+2 \sigma
(1-\sigma)\Gamma(A_1,A_2) + (1-\sigma)^2 \det{A_2}
\eqnl{fi} and
\bqn
\psi(\sigma)= \frac{1}{\det{A_2}}(\sigma^2 \det{A_1}\det{A_2}
+ \sigma(1-\sigma)\tr(A_1 A_2) + (1-\sigma)^2)\,.
\eqnl{psi}
\label{l-shorten}
\el

{\bf Proof of S1.} Recall that the main result in \cite{shorten} claims that the system \eqref{sw-ls} admits a quadratic \lf\
 if and only if $\phi(\sigma)>0$ and $\psi(\sigma)>0$ for every $\sigma\in [0,1]$.
Notice from Lemma~\ref{l-shorten} that $\phi(\sigma)>0$  if and only if either $\Gamma(A_1,A_2)>0$ or the discriminant $\Delta$ of \eqref{fi} is negative. Analogously $\psi(\sigma)>0$ if and only if either $\tr(A_1 A_2)>0$ or the discriminant $\ \tr(A_1 A_2)^2-4\det(A_1)\det(A_2)\ $ of \eqref{psi} is negative. It is therefore clear that the cases considered in {\bf S1} are those satisfying the conditions of \cite{shorten}.
The last statement of {\bf S1} comes from the following series of inequalities
$$\sqrt{\det(A_1)\det(A_2)}\geq \frac12(\tr(A_1)\tr(A_2)-\tr(A_1 A_2))\geq -\frac12\tr(A_1 A_2)\,.$$
This concludes the proof of {\bf S1}.

\vspace{3pt}

{\bf Proof of Corollary~\ref{pos}.} To prove Corollary~\ref{pos} in the case $\det([A_1,A_2])>0$ we use the point \textit{1.(b)} of Lemma~\ref{l-nf-r}. In particular we have that
$\Gamma(A_1,A_2)=\tau_1\tau_2-k>0$ and $\tr(A_1 A_2)=2(\tau_1\tau_2+k)>0$ so that  the conditions of {\bf S1} are satisfied. In the case $\det([A_1,A_2])=0$ the result was already known (see for instance \cite{liberzon-book}), and it can be easily proved by using the normal forms defined in Lemma~\ref{l-nf-r}.

\vspace{3pt}

In what follows we will  always assume $\det([A_1,A_2])<0$.

\vspace{3pt}


{\bf Proof of S2 and S3.} Assume that $\Gamma{(A_1,A_2)} \leq -\sqrt{\det(A_1)\det(A_2)}$. Then a straightforward computation shows that the minimum of $\phi(\sigma)$ is  given by
$$\sigma_0=\frac{\det{A_2} - \Gamma{(A_1,A_2)}}{\det{A_1}+\det{A_2} - 2\Gamma(A_1,A_2)}\in (0,1)\ \mbox{ and }\  \phi(\sigma_0)=\frac{-\Delta}{4( \det A_1+\det A_2-\Gamma{(A_1,A_2))}}\leq 0\,.$$
In particular in the case described by {\bf S2} we have $\phi(\sigma_0)<0$ and therefore the matrix $\sigma_0 A_1+(1-\sigma_0)A_2$ has a positive real eigenvalue, so that the system is unbounded (see Remark \ref{r-conv}).

Similarly when $\Delta=0$ and $\Gamma(A_1,A_2)<0$ we have $\phi(\sigma_0)=0$ so that the system cannot be \guas . In this case to prove that the system is uniformly stable it is possible to show that
the system  admits the following non strict quadratic \lf :
$$V(x)= V(x_1,x_2)=x_1^2 +
\frac{(sign(\delta_{A_1}) sign(\delta_{A_2})-F^2)^2}{4 F^2(\tau_{A_1} F-\tau_{A_2}
sign(\delta_{A_1}))^2}x_2^2  $$

{\bf Proof of S4.}
First observe that, under the conditions of {\bf S4}, we have  $F<-1$ and $k<0$ since, when $A_1$ and $A_2$ are in normal form,  $\tr(A_1 A_2)=F+\frac{sign(\delta_{A_1}\delta_{A_2})}F+2\tau_1\tau_2=2(k+\tau_1\tau_2)<0$ and $|F|>1$.

To prove {\bf S4} we introduce the set of points where the vector fields $A_1x$ and $A_2x$ are parallel:
$$\mathcal{Z}=\{x\in\R^2\,:\,Q(x)=0\}\,,$$
where $Q(x):=\det(A_1x,A_2x)$.
The discriminant of the quadratic function $Q(x)$ coincides with $\Delta$.
Since $\Delta>0$ then $\mathcal{Z}$ consists on a pair of noncoinciding straight
lines passing through the origin. Take a point $x\in\Z\setminus\{0\}$. We say that $\Z$
is {\it direct} (respectively, {\it inverse}) in $x$ if $A_1x$ and $A_2x$ have
the same (respectively, opposite) versus.
We have the following lemma.
\bl
If $\mathcal{Z}$ is direct (resp. inverse) in $x_0\in \Z\setminus\{0\}$ then $\mathcal{Z}$ is direct (resp. inverse) in every point of $\Z\setminus\{0\}$. Moreover in the case {\bf S4} we have that $\Z$ is always direct.\label{dir-inv}
\el
{\bf Proof of Lemma~\ref{dir-inv}.} Let $\Z= D_1\cup
D_2  \ \ \text{where} \ \ D_1,\ D_2$ are straight lines passing through the origin.
Let us observe that, if $x\in D_i$
$$\exists\,\al_i\in\R\quad s.t.\quad A_2A_1^{-1}A_1x=A_2x=\al_i A_1x\,,$$
i.e.  $\al_i$ is an eigenvalue of  $A_2A_1^{-1}$ and $A_1x$ belongs to the eigenspace associated to it.
So $\al_1\al_2=\det(A_2 A_1^{-1})=\frac{\det(A_2)}{\det(A_1)}>0$ which implies that $sign(\al_1)=sign(\al_2)$ i.e. $\Z$ is either direct in every point or inverse in every point. On the other hand it is easy to verify that
$A^{-1}=\frac1{\det{A_1}}(2\tau_1 Id-A_1)$, where $Id$ denotes the identity matrix, which, in the case {\bf S4}, implies
$$\al_1+\al_2=\tr(A_2 A_1^{-1})=\tr\Big(\frac1{\det{A_1}}(2\tau_1 A_2-A_2 A_1)\Big)=\frac{2\Gamma(A_1,A_2)}{\det(A_1)}>0$$
so that $\Z$ is direct. \qed

Let $m_i$ be the slope of $D_i$, for $i=1,2$. Then, if $v_i$ is a vector spanning $D_i$, the orientation of the vector $A_1v_i$ with respect to the radial direction is determined by the quantity
$$sign(\det \left(A_1 v_i,v_i\right))= sign(m_i^2-sign(\delta_{A_1})).$$
Similarly, the orientation of the vector $A_2v_i$ with respect to the radial direction is given by
$$sign(\det \left(A_2 v_i,v_i\right))=sign\left(\frac{m_i^2 sign(\delta_{A_2})-F^2}{F}\right)= sign(F^2-m_i^2 sign(\delta_{A_2})).$$
\bl
If $\Z$ is direct, i.e. if $\Gamma(A_1,A_2)>0$, it must be
$$\ sign(m_i^2-sign(\delta_{A_1}))=sign(F^2-m_i^2 sign(\delta_{A_2}))=+1$$
\label{clock} \el {\bf Proof of Lemma~\ref{clock}.} Since
$\Gamma(A_1,A_2)>0$ by the previous equalities we get that
$\varepsilon:=\ sign(m_i^2-sign(\delta_{A_1}))=sign(F^2-m_i^2
sign(\delta_{A_2}))$. If $\eps=0$ we are in the conditions of {\bf
S1}, since $[A_1,A_2]=0$. If $\eps=-1$ then it must be
$sign(\delta_{A_1})=sign(\delta_{A_2})=1$. In this case we have
$F^2<m_i^2<1$ which is impossible since $|F|>1$.\qed

As a consequence the vectors $A_i x$ point in the clockwise sense
for every $x\in\Z$. This property allows to define the main tool for
checking the stability of \eqref{sw-ls} under the conditions of {\bf
S4}.

\begin{figure}
\begin{center}
\includegraphics[width=0.37\textwidth]{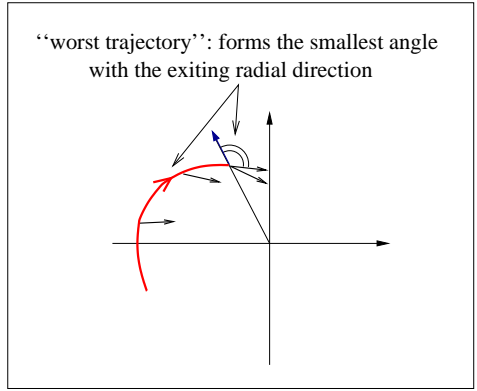}
\caption{The worst trajectory}
\label{worst}
\end{center}
\end{figure}

\bdeff Assume that we are in the conditions of {\bf S4} and under
the normal forms of Lemma~\ref{l-nf-r}. Fix
$x_0\in\R^2\setminus\{0\}$. The worst-trajectory $\g_{x_0}$ is the
trajectory of \eqref{sw-ls}, based at $x_0$, and  having the following
property. At each time $t$, $\dot \g_{x_0}(t)$ forms the smallest
angle in clockwise sense with the exiting radial direction. \edeff
Figure~\ref{worst}  expresses graphically the meaning of the previous definition.

It is clear that the worst trajectory always rotates clockwise
around the origin when $\delta_{A_i}\leq 0$ for some $i\in\{0,1\}$.
If $\delta_{A_i}>0$ for $i=1,2$ then the eigenvectors of $A_1$ are
$(1,1)^T$ and $(1,-1)^T$, while the eigenvectors of $A_2$ are
$(1,F)^T$ and $(1,-F)^T$. In this case it is easy to check that $m_1
m_2<0$ and therefore, from Lemma~\ref{clock}, without loss of
generality we can assume
$$F<m_2<-1<1<m_1<-F.$$
\begin{figure}
\begin{center}
\includegraphics[width=0.37\textwidth]{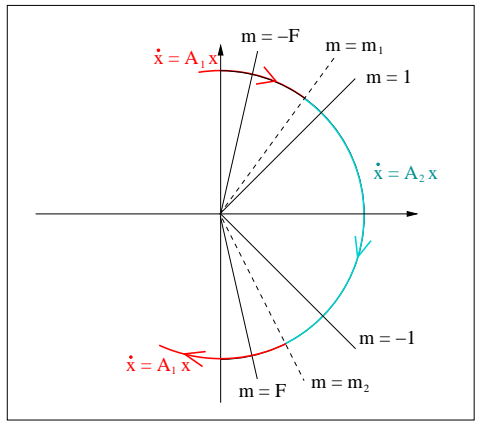}
\caption{Construction of the worst trajectory in the case $\delta_{A_1}>0,\ \delta_{A_2}>0$}
\label{real}
\end{center}
\end{figure}
As a consequence $D_1$ and $D_2$ divide the space into four
connected components, each one intersecting the eigenspace of
exactly one among $A_1$ and $A_2$. This implies that also in this
case the worst trajectory rotates clockwise around the origin. This
trajectory is the concatenation of integral curves of $A_2x$ from
points of $D_1$  to  points of $D_2$ and integral curves of $A_1x$
from points of $D_2$  to  points of $D_1$ (see Figure~\ref{real}).

\begin{figure}
\begin{center}
\includegraphics[width=0.75\textwidth]{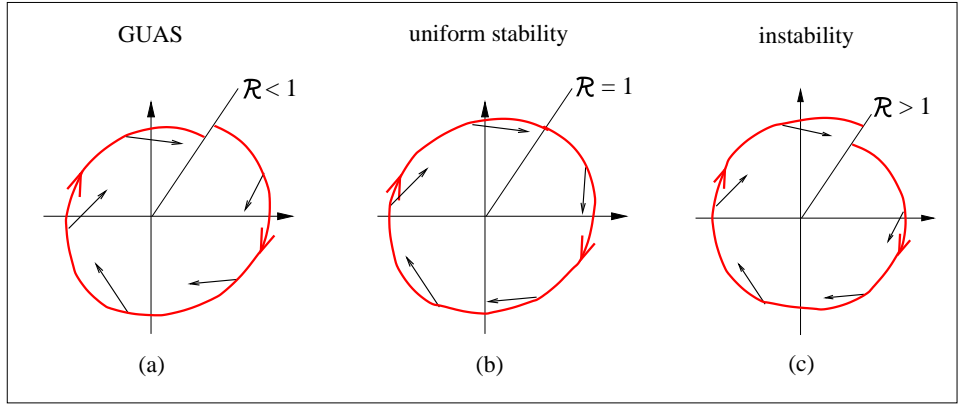}
\caption{The worst trajectory: meaning of the function $\mathcal{R}$}
\label{worst2}
\end{center}
\end{figure}

As explained in the previous papers \cite{sw-balde,sw-1} the
behaviour of the worst trajectory is sufficient to derive the
stability properties of \eqref{sw-ls}. Let us analyse the worst
trajectory $\g_{x_0}(\cdot)$ where $x_0\in D_i$. Assume that $T>0$
is such that $x_1=\g_{x_0}(T)$ is the first intersection point
between the worst trajectory and $D_i$. The worst trajectory tends
to the origin as time goes to infinity if and only if
$\mathcal{R}:=|x_1|/|x_0|<1$, and in this case the system is \guas\ (see Figure~\ref{worst2}~(a)).
It is periodic if and only if $\mathcal{R}=1$, and in this case the
system is uniformly stable but not \guas\ (see Figure~\ref{worst2}~(b)). It blows up if and only if
$\mathcal{R}>1$, and in this case the system is unbounded (see Figure~\ref{worst2}~(c)).

The computation of $\mathcal{R}$ was done in details in
\cite{sw-balde,sw-1}. The formula which is given in
Section~\ref{s-2} is a simpler reformulation of the ones previously
obtained, in terms of our invariants. This concludes the proof of
{\bf S4}.

\section*{Acknowledgements}
Moussa Balde would like to thank the Laboratoire des Signaux et Syst\`emes (LSS - Sup\'elec) for its kind hospitality during the writing of this paper.

{\small
\bibliographystyle{abbrv}
\bibliography{biblio-swit}
}


\end{document}